\documentclass[12pt]{article}
\usepackage{amsmath, amssymb}
\usepackage{amssymb,pb-diagram}
\usepackage{amscd}

\usepackage{fancybox}

\newcommand{\FF}{\mathbb F}

\newcommand{\ZZ}{\mathbb Z}
\newcommand{\PP}{\mathbb P}
\newcommand{\QQ}{\mathbb Q}
\newcommand{\CC}{\mathbb C}

\newcommand{\mcQ}{\mathcal Q}

\newcommand{\tor}{{\mathrm {tor}}}

\newcommand{\red}{{\mathrm {red}}}

\newcommand{\boldc}{\boldsymbol {c}}
\newcommand{\tboldc}{\tilde{\boldsymbol{c}}}
\newcommand{\boldzero}{\boldsymbol {0}}

\newcommand{\MW}{\mathop {\rm MW}\nolimits}

\newcommand{\Div}{\mathop {\rm Div}\nolimits}
\newcommand{\NS}{\mathop {\rm NS}\nolimits}

\newcommand{\Red}{\mathop {\rm Red}\nolimits}
\newcommand{\rank}{\mathop {\rm rank}\nolimits}
\newcommand{\Pic}{\mathop {\rm Pic}\nolimits}

\newcommand{\Supp}{\mathop {\rm Supp}\nolimits}
\newcommand{\Sing}{\mathop {\rm Sing}\nolimits}

\newcommand{\mcB}{\mathop {\mathcal B}\nolimits}

\newcommand{\mcO}{\mathop {\mathcal O}\nolimits}

\newtheorem{thm}{Theorem}[section]
\newtheorem{cor}{Corollary}[section]
\newtheorem{prop}{Proposition}[section]
\newtheorem{lem}{Lemma}[section]
\newtheorem{defin}{Definition}[section]
\newtheorem{exmple}{Example}[section]
\newtheorem{rem}{Remark}[section]
\newtheorem{qz}{Question}[section]

\newcommand{\I}{\mathop {\rm I}\nolimits}

\newenvironment{example}{\begin{exmple}\rm }{\end{exmple}}

\newcommand{\qed}{\hfill $\Box$}

\newcommand{\proof}{\noindent{\textsl {Proof}.}\hskip 3pt}
\newcommand{\proofend}{\qed \par\smallskip\noindent}

\renewcommand{\thesubparagraph}{\theparagraph.\@arabic\c@subparagraph}
\begin{document}
  
\bigskip

\begin{center}

{\Large \bf  On the Abel-Jacobi map of an elliptic surface\\ and  \\ the topology of cubic-line arrangements}

\bigskip

Shinzo BANNAI and 
Hiro-o TOKUNAGA

\end{center}
\abstract{Let $\varphi : S \to C$ be an elliptic surface over a smooth curve $C$ with a section $O$. We
denote its generic fiber by $E_S$. 
For a divisor $D$ on $S$, we canonically associate a $\CC(C)$-rational point $P_D$. In this note,
we give a description of $P_D$ of $E_S$, when the rank of the group of $\CC(C)$-rational points is one.
We apply our description to refine our result on a Zariski pair
for a cubic-line arrangement.}

\section*{Introduction}

Let $E$ be an elliptic curve defined over a field $K$  isomorphic to either $\CC(t)$ or $\CC(t_1, t_2)$ (the rational function field of one variable or two). We
denote their group of $K$-rational points by $E(K)$. Since $E$ can be 
considered the generic fiber of an elliptic surface or $3$-fold, 
they both  have arithmetic and geometric aspects.
In \cite{tokunaga98},  the second author considered the case when $E$ is the generic fiber of a certain elliptic K3 surface, 
and made use of 
$3$-torsions of $E(\CC(t))$ in order to construct Zariski pairs for irreducible sextic curves. Also, in our previous works \cite{bannai-tokunaga15, bannai-tokunaga17,
tokunaga14}, we investigated the case when $E$ is the generic fiber of  certain rational  elliptic surfaces,  and costructed Zariski pairs ($N$-plets) 
 for reducible curves by using non-torsion elements in $E(\CC(t))$. 
 For the case of $\CC(t_1, t_2)$, Cogolludo Agsutin, Kloosterman and Libgober have recently investigated
 $E(\CC(t_1, t_2))$ in order to study toric decomposition of plane curves, which is related to the embedded topology of plane curves
(\cite{cogo-kloosterman, cogo-lib, kloosterman2013, kloosterman2014}).    These results shows that the study of arithmetic aspects of elliptic curves 
over the rational function fields is one of the important tools to study the topology of plane curves. In this article, we continue to study 
topology of plane curves along this line with more emphasis on the arithmetaric aspects, especially the Abel-Jacobi map on an elliptic surface, which we explain below.

Let $\varphi: S \to C$ be an elliptic surface over a smooth projective curve. Throughout this paper, we 
assume that 
(i) $\varphi$ is relatively minimal,
(ii) there exists a section $O: C \to S$, and 
(iii) there exists at least one degenerate fiber.
As for a section $s : C \to S$, we identify $s$ and 
its image, i.e, an irreducible curve meeting any fiber at one point.  We denote the set of sections of $\varphi : S \to C$ by $\MW(S)$.
Note that $\MW(S) \neq \emptyset$ as $O \in \MW(S)$.
 Let $E_S$, $\CC(C)$ and $E_S(\CC(C))$ denote the generic fiber of $S$, the rational function
field of $C$  and the set of $\CC(C)$-rational points of $E_S$, respectively.  
Under our assumption, by restricting a section to the generic fiber, $\MW(S)$  can be canonically identified
with $E_S(\CC(C))$. We identify $O$ with the corresponding rational point. Thus $(E_S, O)$ is an elliptic
curve defined over $\CC(C)$. 

Let $D$ be a divisor on $S$. By restricting $D$ to $E_S$, we have a divisor ${\frak d}$ on $E_S$ defined over 
$\CC(C)$.  By applying Abel's theorem to ${\frak d}$, we have $P_D \in E_S(\CC(C))$ and thus $s(D) \in \MW(S)$ (see
\cite[Lemma~5.1, \S 5]{shioda90} for the explicit description of $s(D)$). 

In our previous articles \cite{bannai-tokunaga15, bannai-tokunaga17,  tokunaga04, tokunaga14}, we studied the properties
of $P_D$ in $E_S(\CC(C))$ such as $p$-divisibility, for odd primes $p$, in order to study the topology of reducible plane curves
with irreducible components of low degrees.
 In this article, we  consider $n$-divisibility of $P_D$ in the case when $\rank E_S(\CC(C))$ $= 1$, i.e., $E_S(\CC(C)) = \ZZ P_o \oplus E_S(\CC(C))_{\tor}$ for some $P_o \in E_S(\CC(C))$.

 \begin{prop}\label{prop:main} {Let $\phi_o : \NS(S) \to \NS(S)\otimes \QQ$ and $\phi : E_S(\CC(C)) \to  \NS(S)\otimes \QQ$
   be the homomorphism defined in \S 2. Suppose that $\rank E_S(\CC(C)) = 1$.  Let $n$ be
 an integer such that $P_D = nP_o + P_{\tau}$, $P_{\tau} \in E_S(\CC(C))_{\tor}$. Then we have
\[
 n^2  =  - \frac {\phi_o(D)\cdot \phi_o(D)}{\langle P_o, P_o\rangle}, \quad
 n  =  -\frac {\phi_o(D)\cdot \phi(P_o)}{\langle P_o, P_o \rangle}
\]
 where $\cdot$ and $\langle \, , \, \rangle$ mean the intersection and height pairing, respectively.
 }
 \end{prop}

\begin{rem}{\rm 
As we see in Lemma~\ref{lem:fundamental}, \S 2, $\phi_o(D) = \phi(P_D)$. Hence $\phi_o(D)$ can be considered as (almost) a ^^ class' in the Picard group of the generic fiber.
}
\end{rem}

An explicit form for the right hand side in the above formula are given in \S 2.  In \S 3 we develop a method to determine the 
contribution from the torsion part  from $\phi_o(D)$.  Hence for the rank one case, it is possible to describe $P_D$ completely. 

\medskip

In the remaining part of this article, we consider an application of Proposition~\ref{prop:main} in the investigation of  the embedded topology of reducible
 plane curves.  As we have seen in our previous papers, properties of $P_D$ in $E_S(\CC)$  played important
 role in order to study the topology of plane curves which arise from $D$. 
 In this article, we compute $n$-divisibility of $P_D$ for a certain trisection  and apply it to
   refine our result for a Zariski pair  given in \cite{bty17} as follows:

 Let $(\mcB^1, \mcB^2)$ be the Zariski pair for a nodal cubic $E$ and  four lines considered in \cite{bty17}, i.e., the one with
 Combinatorics 1-(b).  Namely, it is as follows:
 \medskip

 {\bf Combinatorics 1-(b).} Let $E$, and $L_i$ ($i = 0,  1, 2, 3$) be as below and
we put $\mcB = E  +  \sum_{i=0}^3 L_i$:

\begin{enumerate}

\item[(i)] $E$:  a  nodal cubic curve.

\item[(ii)] $L_0$: a transversal line to $E$ and we put $E \cap L_0 = \{p_1, p_2, p_3\}$.

\item[(iii)] $L_i$: a line through $p_i$ and tangent to  $E$ at a point $q_i$ distinct from  $p_i$
$(i =1, 2, 3)$.

\item[(iv)] $L_1, L_2$ and $L_3$ are not concurrent.

\end{enumerate}

By taking the group structure of $E\setminus \{the \, \, node\}$ into account, we infer that $q_i$ $(i = 1, 2, 3)$ are either collinear or not. 
For $\mcB$ with Combinatorics 1-(b), we call it Type I (resp. Type II) if $q_1, q_2$ and $q_3$ are collinear (resp. not collinear)

 Then 
  with terminologies and notation for $D_{2n}$-covers given in \cite{act, bannai-tokunaga15}, we have
 
 \begin{thm}\label{thm:refine}{Let $n$ be an integer $\ge 3$.
 \begin{enumerate}
 
  \item[(i)] If $\mcB$ is of Type I, there exists a $D_{2n}$-cover $\pi : X_n \to \PP^2$ branched at $2(\sum_{i= 0}^3 L_i) + n E$ for any $n$.
 
 \item[(ii)]  If $\mcB$ is of Type II, there exists a $D_{2n}$-cover $\pi : X_n \to \PP^2$ branched at $2(\sum_{i= 0}^3L_i) + n E$ for $n = 4$ only.

 \end{enumerate}
 }
 \end{thm}
 
 \begin{cor}\label{cor:refine}{Let $(\mcB^1, \mcB^2)$ be  a pair of plane curves with Combinatorics 1-(b) such that their Types are distinct. Then
  both of the fundamental groups $\pi_1(\PP^2\setminus \mcB^j, \ast)$ ($j = 1, 2$) are non-abelian and
 there exist no homeomorphisms between $(\PP^2, \mcB^1)$ and $(\PP^2, \mcB^2)$.
 }
 \end{cor}
 
 \begin{rem}{\rm For $\mcB$ of Type I, we denote the line through $q_i$ $(i = 1, 2, 3)$ by $\overline {L}$. Then we infer that $\sum_{i=1}^3L_i $ is a member
 of the pencil generated by $E$ and $L_0 + 2\overline{L}$.  This shows that there exists a $D_{2n}$-cover of $\PP^2$ branched at $2(E + L_0) + n(\sum_{i=1}^3L_i)$
 for any $n\,  (\ge 3)$. If $p_i$ is not an inflection point of $E$,  there exist just two lines through $p_i$ such that they are tangent to
 $E$ at different points from $p_i$. Hence  we infer that for the pair $(\mcB^1, \mcB^2)$ for
 Combinatorics 1-(b) given in
 \cite{bty17}, $\mcB^1$ (resp. $\mcB^2$) is Type II (resp. Type I).
 }
 \end{rem}

\begin{rem}{ \rm For the explicit example for $(\mcB^1, \mcB^2)$ given in \cite{bty17}, the non-abelianness for $\pi_1(\PP^2\setminus \mcB^j, \ast)$
$(j = 1, 2)$ was first pointed out 
by E. Artal Bartolo. In this note, we prove that the same is true for any curve with the same combinatorics.
}
\end{rem}

\bigskip

 
\section{Preliminaries on elliptic surfaces}\label{sec:junbi}

We refer to \cite{kodaira}, \cite{miranda-basic}, \cite{miranda-persson} and 
\cite{shioda90} for details.
In this article, an {\it elliptic surface} always means the one introduced in the Introduction.
 We denote a subset of $C$ over which $\varphi$ has degenerate fibers by $\Sing(\varphi)$.
 $\Red(\varphi)$ means a  subset of $\Sing(\varphi)$ consisting of a point 
 $v \in \Sing(\varphi)$ such that $\varphi^{-1}(v)$ is reducible.
  For $v \in \Sing(\varphi)$, we denote the corresponding fiber by $F_v = \varphi^{-1}(v)$. 
 The irreducible decomposition of $F_v$ is denoted by 
 \[
 F_v = \Theta_{v, 0} + \sum_{i=1}^{m_v-1} a_{v,i}\Theta_{v,i}
 \]
 where $m_v$ is the number of irreducible components of $F_v$ and $\Theta_{v,0}$ is the
 irreducible component with $\Theta_{v,0}O = 1$. We call $\Theta_{v,0}$ the {\it identity
 component}.  In order to describe types of singular fibers, we use
 Kodaira's symbol.  We label irreducible components  of singular fibers as in \cite[p.81-82]{tokunaga12}.

 Let $\MW(S)$ be the set of sections of $\varphi : S \to C$. By our assumption,
 $\MW(S) \neq \emptyset$, as $O \in \MW(S)$. By regarding $O$ as the zero 
 element,  $\MW(S)$ is equipped with the structure of an abelian group through fiberwise addition.

  Let $E_{S}$ be the generic fiber of 
 $\varphi : S \to C$. We can regard $E_{S}$ as a curve of genus $1$ over
 $\CC(C)$, the rational function field of $C$ and that under our setting, as $S$ is the 
 the Kodaira-N\'eron model of $E_S$,  $\MW(S)$ is identified with
 the set of $\CC(C)$-rational points, $E_S(\CC(C))$, of $E_{S}$.
 
 Let $\NS(S)$ be the N\'eron-Severi group of $S$. Under our assumption on $S$, $\NS(S)$ is
 torsion free by \cite[Theorem~1.2]{shioda90}.  Let $T_{\varphi}$ be the subgroup of 
 $\NS(S)$ generated by $O$, a fiber $F$ of $\varphi$ and $\Theta_{v,i}$ $(v \in
 \Red(\varphi),  1 \le i \le m_v-1)$. In \cite[Theorem~1.3]{shioda90}, by taking Abel's theorem on 
 $E_{S}$ into account, an isomorphism $\overline{\psi} : \NS(S)/T_{\varphi} 
 \to E_S(\CC(C))$ of  abelian groups is given as follows:
 
 We first define a homomorphism $\psi$ from the group of divisors ${\mathrm {Div}}(S)$ to 
 $\Pic_{\CC(C)}^0(E_S) \cong E_S(\CC(C))$ by 
 \[
 \psi : {\mathrm {Div}}(S) \ni D \mapsto \alpha(D|_{E_{S}} - (DF)O|_{E_{S}}) \sim_{E_S} P_D  - O,
 \]
where $\alpha$ is the Abel-Jacobi map on $E_S$ and $\sim_{E_S}$ denotes the linear 
equivalence on $E_S$. By \cite[Lemma~5.2]{shioda90}, 
$\psi$ induces a group isomorphism (\cite[Theorem~1.3]{shioda90})
\[
\overline{\psi} : \NS(S)/T_{\varphi} \to E_S(\CC(C)).
\] 
 We denote the  section corresponding to $P_D$ by  $s(D)$. By \cite[Lemma~5.1]{shioda90}, we have a relation in 
 $\NS(S)$: 
 \[
(\ast) \quad D \approx s(D) + (d-1)O + nF + \sum_{v\in \Red(\varphi)}\sum_{i=1}^{m_v-1}b_{v,i}\Theta_{v,i},
\]
where $\approx$ denotes the algebraic equivalence between divisors, and $d, n$ and $b_{v,i}$
are integers defined as follows:
\[
d = D\cdot F \qquad n = (d-1)\chi({\mathcal O}_S) + O\cdot D - s(D)\cdot O,
\]
and
\[
\left [ \begin{array}{c}
          b_{v,1} \\
          \vdots \\
          b_{v, m_v-1} \end{array} \right ] = A_v^{-1}\left (\boldc(v, D) - \boldc(v, s(D)) \right ), 
\]
 where, for $D \in \mathrm{Div}(S)$,  we put    
\[          
 \boldc(v, D) :=\left [
 \begin{array}{c}
 D\cdot\Theta_{v,1}\\
 \vdots \\
   D\cdot \Theta_{v,m_v-1}
 \end{array} \right ],
\]
and $A_v$ is the intersection matrix $(\Theta_{v,i}\Theta_{v, j})_{1\le i, j \le m_v-1}$.

\begin{rem}\label{rem:restriction}{\rm
\begin{enumerate}

\item[(i)] If $D$ is a section (i.e, $D \in \MW(S)$), the above relation $(\ast)$ becomes trivial.

\item[(ii)] Entries of  $A_v^{-1}$ are not necessarily integers. On the other hand,
the relation $(\ast)$ is a relation between two divisors of $\ZZ$-coefficients. This impose
some restriction on $D$ and $s(D)$ at which irreducible components of $F_v$, $D$ and $s(D)$
meet.  One of useful facts is a lemma below.

\end{enumerate}
}
\end{rem}
 
 \begin{lem}\label{lem:key-1}{ If $A_v^{-1}\boldc(v, D) \in \ZZ^{\oplus m_v - 1}$, then $\boldc(v, s(D)) = \boldzero$.
 }
 \end{lem}
 
 \proof  Since $s(D)\cdot F = 1$, $\tboldc(v, s(D))$ has a unique entry with $1$ and other entries are $0$. Also,
 $s(D)$ meets the component $\Theta_{v, i}$ with $a_{v, i} = 1$ or $\Theta_{v,0}$. Hence if $\boldc(v, s(D)) \neq \boldzero$,
 $A_v^{-1}\left (\boldc(v, D) - \boldc(v, s(D)) \right ) \not\in \ZZ^{\oplus m_v - 1}$. \qed 
 
\section{Proof of  Proposition~\ref{prop:main}}\label{sec:prop-proof}
 
 Put $\NS_{\QQ}:= \NS(S)\otimes \QQ$ and $T_{\varphi, \QQ} := T_{\varphi}\otimes \QQ$.
As $\NS(S)$ is torsion free by \cite[Theorem~1.2]{shioda90}, there is no big difficulty in considering $\NS_{\QQ}$.
 By using
 the intersection pairing, we have an orthogonal decomposition
 \[
 \NS(S)_{\QQ} = T_{\varphi, \QQ} \oplus T^{\perp}_{\varphi, \QQ}.
 \]
 Let $\phi$ denote the homomorphism from $E_S(\CC(C))$ to $\NS(S)_{\QQ}$ given in 
 \cite[Lemma~8.2]{shioda90}. Also we define a homomorphism $\phi_o$ from $\Div(S)$ to $T_{\varphi, \QQ}^{\perp} ( \subset \NS(S)_{\QQ})$
 by the composition:
 \[
 \phi_o: \Div(S) \to \NS(S) \to T_{\varphi, \QQ}^{\perp} \subset \NS(S)_{\QQ} ,
 \]
 the last morphism is the projection.  Explicitly, for $D \in \Div(S)$, $\phi_o$ is given by
 \[
 (\ast\ast) \quad  \phi_o(D) = D - dO - (d\chi + (O\cdot D))F - \sum_{v \in \Red(\varphi)}\FF_v A_v^{-1}\boldc(v, D), 
 \]
  where $d = D\cdot F, \chi = \chi({\mathcal O}_S)$ and $\FF_v = [\Theta_{v, 1}, \ldots,  \Theta_{v, m_v-1} ]$.
 Here we have the following lemma on $\phi$ and $\phi_o$:
 
 \begin{lem}\label{lem:fundamental}
 \begin{enumerate}
 
 \item[(i)] For $P \in E_S(\CC(C))$ and its corresponding section $s_P$, we have $\phi(P) = \phi_o(s_P)$.
 
 \item[(ii)] For $D \in \Div(S)$ and its corresponding point $P_D \in E_S(\CC(C))$,  we have
 $\phi(P_D) = \phi_o(D)$.
 
 \end{enumerate}
 \end{lem}
 
 \proof The statement (i) follows from definition. For (ii), our statement follows from the relation $(\ast)$ in the previous
 subsection. \qed
 
 \bigskip

 In \cite{shioda90}, a $\QQ$-valued bilinear form $\langle \, , \, \rangle$  called the {\it height paring} on $E_S(\CC(C))$ is defined by  
 $\langle P_1, P_2 \rangle := -\phi(P_1)\cdot\phi(P_2)$ (see  \cite{shioda90} for details). 
 By Lemma~\ref{lem:fundamental} (ii),  we have
\[
\langle P_D, P_D \rangle = - \phi_o(D)\cdot \phi_o(D), \quad \langle P_D, P_o\rangle  = - \phi_o(D)\cdot\phi(P_o).
\]
As $\langle P_D, P_D \rangle = n^2\langle P_o, P_o \rangle$ and $\langle P_D, P_o \rangle = n\langle P_o, P_o\rangle$, 
we have our statement in Proposition~\ref{prop:main}. Also by computing
the intersection pairing explicitly, we have
\begin{eqnarray*}
 \phi_o(D)\cdot \phi_o(D) & = & D^2 -2d D\cdot O - d^2\chi - \sum_{v \in \Red(\varphi)} {}^t\!\boldc(v, D)A_v^{-1}\boldc(v, D) \\
 \phi_o(D)\cdot \phi(P_o) & = & (D -d O)\cdot s_{P_o} - d\chi  -O\cdot D - \sum_{v \in \Red(\varphi)} {}^t\!\boldc(v, s_{P_o})A_v^{-1}\boldc(v, D) \\
 \langle P_o, P_o \rangle & =& 2\chi + 2s_{P_o}\cdot O  + \sum_{v \in \Red(\varphi)} {}^t\!\boldc(v, s_{P_o})A_v^{-1}\boldc(v, s_{P_o}).
 \end{eqnarray*} 
 
 Note that we do not need any data of $P_D$ in order to compute $\phi_o(D)\cdot \phi_o(D)$.

 
 
\section{The torsion part of $P_D$}\label{sec:torsion}

\subsection{The homomorphism $\gamma_{\NS}$}\label{sec:section2}

For a reducible singular fiber $F_v = \sum_i a_{v, i}\Theta_{v,i}$ $(v \in \Red(\varphi))$,
we denote a subgroup generated by $\Theta_{v,1}, \ldots, \Theta_{v, m_v-1}$ by $R_v$.
Let $R_v^{\vee}$ be the dual of $R_v$, which can be embedded into $R_v\otimes \QQ$
by the intersection pairing. Under this circumstance, $R_v^{\vee}$ can be regarded as
a subgroup generated by the columns of $A_v^{-1}$.

\begin{defin}\label{def:gammaNS}{\rm 
 We define a map $\gamma_{\NS}$ from
$\NS(S)$ to $\oplus_{v \in \Red(\varphi)} R_v^{\vee}$ by
\[
\gamma_{\NS} : \NS(S) \ni D \mapsto \left ( -A_v^{-1} \boldc(v, D)\right )_{v \in \Red(\varphi)}
\in \bigoplus_{v \in \Red(\varphi)}R_v^{\vee},
\]
where $\boldc(v, D)$ as in  \S~\ref{sec:junbi}. 
We denote the induced
 map from $\NS(S)$ to $\oplus_{v\in \Red(\varphi)}R_v^{\vee}/R_v$ by
 \[
 \overline{\gamma}_{\NS} : \NS(S) \to \oplus_{v\in \Red(\varphi)}R_v^{\vee}/R_v.
 \]

}
\end{defin}

\begin{lem}\label{lem:gammaNS}{Both $\gamma_{\NS}$ and 
$\overline{\gamma}_{\NS}$ are group homomorphisms.
}
\end{lem}

\proof Since $\Theta_{v,i}\cdot (aD_1 + bD_2) = a\Theta_{v,i}\cdot D_1 + b\Theta_{v,i}\cdot D_2$ 
$(D_1, D_2 \in \NS(S), a, b \in \ZZ)$, our statement is immediate. \proofend 

\begin{lem}\label{lem:basic}{ Let $\psi : \Div(S) \to E_S(\CC(C))$ be the homomorphism in the 
previous section. For $D_1, \, D_2 \in \Div(S)$, if $\psi(D_1) = \psi(D_2)$, then 
$\overline{\gamma}_{\NS}(D_1) = \overline{\gamma}_{\NS}(D_2)$, where we identify $D_i$ with
its algebraic equivalence class.
}
\end{lem}

\proof Put $s(D_i)$ be the corresponding sections to $\psi(D_i)$ $(i = 1, 2)$. Then by $(\ast)$ in \S~\ref{sec:junbi}, we have
\begin{eqnarray*}
D_i & \approx & s(D_i) + (d_i-1)O + n_i F +  \sum_{v\in \Red(\varphi)}\FF_v(-A_v)^{-1}(\boldc(v, D_i) - 
  \boldc(v, s(D_i))),
\end{eqnarray*}
where $d_i = D_iF$, $n_i = (d_i -1)\chi(\mcO_X) + O\cdot D_i - O\cdot s(D_i)$. Since $s(D_1) = s(D_2)$,
we have 
\begin{eqnarray*}
D_1 - D_2 & \approx &  (d_1-d_2)O + (n_1 - n_2)F +  \\
& &  \sum_{v\in \Red(\varphi)}
\FF_v(-A_v)^{-1}(\boldc(v, D_1) - \boldc(v, s(D_2))).
\end{eqnarray*}
In the above equivalence, all coefficients of $O$, $F$, and $\Theta_{v, j}$'s are integers.
Hence all the entries of $(-A_v)^{-1}(\boldc(v, D_1) - \boldc(v, D_2))$ 
$(v \in \Red(\varphi))$ are integers. Since $R_v^{\vee}$ can be regarded as a $\ZZ$-module
obtained by adding column vectors of $(-A)_v^{-1}$, we infer that 
 $(-A_v)^{-1}(\boldc(v, D_1) - \boldc(v, D_2)) \in R_v$ for $\forall v \in \Red(\varphi)$. Hence
 we have 
 \[
 \overline{\gamma}_{\NS}(D_1) = \overline{\gamma}_{\NS}(D_2).
 \]
  \proofend
  
  \begin{rem}{\rm For a singular fiber $F_v = \sum_i a_{v,i}\Theta_{v,i}$, we  put
  \[
  F_v^{\sharp} = \cup_{a_{v,i} = 1}\Theta_{v,i}^{\sharp},
  \]
  where $\Theta_{v,i}^{\sharp} := \Theta_{v,i} \setminus (\mbox{singular points of 
  $(F_v)_{\red}$})$. By \cite[\S 9]{kodaira}, $F_v^{\sharp}$ has a structure of 
  an abelian group, and we define an finite abelian group $G_{F_v^{\sharp}}$ as in
  \cite[p.81-82]{tokunaga12}.
Roughly speaking, $G_{F_v^{\sharp}}$ is a group given by the indices of the
irreducible components of $F_v$. Put 
$G_{\Sing(\varphi)}:= \sum_{v\in \Sing{\varphi}} G_{F_v^{\sharp}}$ and we 
define a homomorphism $\gamma : \MW(S) \to G_{\Sing(\varphi)}$, which
describes at which irreducible component each section section meets. The
homomorphism $\overline{\gamma}_{\NS} : \NS(S) \to 
\oplus_{v \in \Red(\varphi)} R_v^{\vee}/R_v$ can be considered as a generalization of $\gamma$. In fact, $R_v^{\vee}/R_v$ is canonically isomorphic
to $G_{F_v^{\sharp}}$ and $\overline{\gamma}_{\NS}$ and $\gamma$ coincide
for sections.
}
\end{rem}

\begin{example}

To illustrate the morphism $\overline{\gamma}_{\NS}$ and the isomorphism between $R_v^\vee/R_v$ and $G_{F_v^\sharp}$ in more detail, let us look in to the case where $S$ has a unique reducible singular fiber $F_v$  of type $\I_0^\ast$. We will relabel the fiber components of $\I_0^*$ as
\[
 \Theta_{v,0}, \Theta_{v,1}=\Theta_{v,01}, \Theta_{v,2}=\Theta_{v,10}, \Theta_{v,3}=\Theta_{v,11}, \Theta_{v,4}.
\]
In this case, we have
\[
A_v=(\Theta_{v,i}\Theta_{v,j})=\left[\begin{array}{cccc} -2 & 0 & 0 & 1\\ 0& -2 & 0 & 1 \\ 0 & 0 & -2 & 1 \\ 1 & 1 & 1 & -2\end{array}\right], A_v^{-1}=
\left[ \begin {array}{cccc} -1&-1/2&-1/2&-1\\ -1/2&
-1&-1/2&-1\\ -1/2&-1/2&-1&-1\\ -1&-1&-1&-2\end {array} \right]. 
\]
Let ${\bf x}={}^t\left[\begin{array}{cccc} x_1 & x_2 & x_3 & x_4\end{array}\right]$ and suppose that 
\[
-A_v^{-1}{\bf x}=-\left[ \begin {array}{cccc} -1&-1/2&-1/2&-1\\ -1/2&
-1&-1/2&-1\\ -1/2&-1/2&-1&-1\\ -1&-1&-1&-2\end {array} \right]\left[\begin{array}{c} x_1 \\ x_2 \\ x_3 \\ x_4\end{array}\right]= \left[\begin{array}{c} a \\ b \\ c \\ d\end{array}\right]\in R_v
\]
for some integers $a, b, c, d\in \mathbb{Z}$.
Then, we have 
$
x_1=2a-d$,
$x_2=2b-d$,
$x_3=2c-d$,
$x_4=-a-b-c+2d$, which implies that $x_1, x_2, x_3$ must have the same parity. Conversely, if $x_1, x_2, x_3$ have the same parity, $-A_v^{-1}{\bf x}\in R_v$. From these facts, we see that in $R_v^\vee/R_v$, any  $-A_v^{-1}{\bf x}$ is equivalent to one of ${\bf 0}$, $-A_v^{-1}{\bf e}_1$, $-A_v^{-1}{\bf e}_2$, $-A_v^{-1}{\bf e}_3$, where ${\bf e}_1, {\bf e}_2, {\bf e}_3$ are the first three of the standard basis vectors. Also we have $(-A_v)^{-1}{\bf e}_i+(-A_v)^{-1}{\bf e}_j=(-A_v)^{-1}{\bf e}_k$, ($\{i,j,k\}=\{1,2,3\}$) and $2(-A_v)^{-1}{\bf e}_i={\bf 0}, (i=1,2,3)$. Therefore, we have $R_v^\vee/R_v\cong (\ZZ/2\ZZ)^{\oplus 2}\cong G_{F_v^\sharp}$.
\end{example}

\subsection{The determination of the torsion part}\label{sec:torsion}
As noted above, Proposition \ref{prop:main} allows us to compute the coefficient of the $P_o$ part of $P_D$. In this subsection we study the homomorphism $\gamma_{\NS}$ in order to determine the torsion part $P_\tau$ of $P_D$. 

Let $P_\tau, P_{\tau^\prime}\in E_S(\CC(C))$ be torsion points and let $s_\tau, s_{\tau^\prime}$  be the corresponding sections. Then we have:

\begin{lem}\label{lem:torsion}
$\bar{\gamma}_{\NS}(s_\tau)= \bar{\gamma}_{\NS}(s_{\tau^\prime})\Leftrightarrow P_\tau=P_{\tau^\prime}$
\end{lem}

\proof
We first recall that $\gamma$ and $\bar{\gamma}_{\NS}$ coincide for sections.
Suppose $\bar{\gamma}_{\NS}(s_\tau)= \bar{\gamma}_{\NS}(s_{\tau^\prime})$. Then $s_\tau\approx s_{\tau^\prime}$ by ($\ast$), which implies that $P_\tau\sim_{E_S} P_{\tau^\prime}$ as divisors on $E_S$. Hence we have $P_\tau=P_{\tau^\prime}$. The converse is obvious from the definition of $\gamma_{\NS}$. \qed

From this Lemma~\ref{lem:torsion},  it is enough to determine $\bar{\gamma}_{\NS}(s_\tau)$ to determine $P_\tau$, but since $\gamma$  is a homeomorphism, we have
\[
\bar{\gamma}_{\NS}(s_\tau)=\bar{\gamma}_{\NS}(P_D)-\bar{\gamma}_{\NS}(nP_o)
\]
which enables us to compute $\bar{\gamma}_{\NS}(s_\tau)$.


\section{A rational elliptic surface attached to general $4$ lines}\label{sec:res-4lines}

In this section, we apply the above discussions to compute $P_D$ for certain $D$ in the case of a rational elliptic surfaces that is associated to an arrangement of four non-concurrent lines $L_i$, $i=0,1,2,3$. The computations will be used to prove Theorem \ref{thm:refine} in the next section.

Let $z_o$ be a general point of $L_0$, $Q=L_0+L_1+L_2+L_3$ and consider the rational elliptic surface $S_{Q,z_o}$ associated to $Q$ and $z_o$. (For the details of the construction of $S_{Q,z_o}$, see \cite[II]{bannai-tokunaga15}.) Let $L_0\cap L_i=p_i$ ($i=1,2,3$) and $L_i\cap L_j=p_{ij}$ ($i\not=j$). By the construction, $S_{Q, z_o}$  is a rational elliptic surface whose set of  reducible singular fibers is of type $\I_0^\ast$, $3\I_2$. The $\I_0^\ast$ fiber arises from  the preimage of the line $L_0$ and the $\I_2$ fibers arise from the preimages of the lines through $z_o$ and $q_i$. We denote the  $\I_0^\ast$ fiber by $F_{v_0}$  and its components by.
\[
 \Theta_{\infty,0}, \Theta_{\infty,1}=\Theta_{\infty,01}, \Theta_{\infty,2}=\Theta_{\infty,10}, \Theta_{\infty,3}=\Theta_{\infty,11}, \Theta_{\infty,4}
\]
as in Section \ref{sec:torsion}, and for $i=1,2,3$,  we label the $\I_2$ fiber corresponding to the line $\overline{z_oq_i}$ by $F_i$ and its components by $\Theta_{i,0}, \Theta_{i,1}$. Also, by \cite{oguiso-shioda}, $\MW(S _{Q,z_o})\cong A_1^\ast\oplus(\ZZ/2\ZZ)^{\oplus 2} $.

Next we consider configurations consisting of an irreducible cubic $E$ and $Q=L_0+L_1+L_2+L_3$ satisfying the following:
\begin{itemize}
\item $E$ passes through the three points $p_1, p_2, p_3$.
\item For $i=1,2,3$, $E$ is tangent to $L_i$ at a point $q_i$ distinct from $p_i, p_{ij}$.
\end{itemize} 

Note that if $E$ is a nodal cubic, $E + \sum_{i= 0}^3L_i$ has Combinatorics 1-(b). Suppose $E$  is a splitting curve with respect to $Q$ and let $E^\pm$  be the irreducible components of the strict transform of $E$ in $S_{Q, z_o}$. 
Our goal is to compute $P_{E^\pm}$. 
We first recall  the diagram 
 \[
\begin{CD}
S'_{Q} @<{\mu}<< S_{Q} @<{\nu_{z_o}}<<S_{Q, z_o} \\
@V{f'_{Q}}VV                 @VV{f_{\mcQ}}V         @VV{f_{Q, z_o}}V \\
\PP^2@<<{q}< \widehat{\PP^2} @<<{q_{z_o}}< (\widehat{\PP^2})_{z_o},
\end{CD}
\]
that appears in the construction of $S_{\mcQ, z_o}$ given in \cite[Introduction]{bannai-tokunaga17}. Note that the covering transformation of $f_{Q, z_o}$ induces
the inversion morphism on  the generic fiber $E_{S_{Q, z_o}}$ and we have $P_{E^-} = -P_{E^+}$.  
Let $P_{E^+}=nP_o+P_\tau$, where $P_o$ is a generator of the $A_1^\ast$ part of $E_S(\CC(t))$.  We may assume $n \ge 0$ after relabeling $E^{\pm}$, suitably.
First, from the data of the intersection of $E$ and $L_i$ ($i=0,1,2,3$) we have  
\[
\boldc(v_0, E^+)={}^t\left[\begin{array}{cccc} 1 & 1 & 1 & 0\end{array}\right], \boldc(v_i, E^+)=0\,(i=1,2,3).
\] 
where $v_0$ is the $\I_0^\ast$ fiber and $v_i$ are the $\I_2$ fibers. 
Now, from Proposition \ref{prop:main}, we have
\[
n^2=-\dfrac{\phi_o(E^+)\cdot\phi_o(E^+)}{\langle P_o, P_o \rangle}=-2((E^+)^2-3).
\]

  Since $E$ is a cubic and passes through $p_1, p_2, p_3$, the strict transform of $E$ in $ (\widehat{\PP^2})_{z_o}$ has self-intersection number 6. Hence, we have 
\[
(E^++E^-)^2=2\cdot 6=12
\]
by which we obtain 
\[
(E^+)^2=(E^-)^2=6-E^+\cdot E^-\,.
\]
Now, if $E$ is smooth, then $E^+\cdot E^-=3$ and we have $(E^+)^2=3$, which implies $n=0$.

Next, if $E$ is a nodal cubic, there are two possibilities for $E^+\cdot E^-$, namely $E^+\cdot E^-=3$ or $5$, depending on the data of the preimage of the node. If the preimage of the node becomes nodes of $E^+, E^-$, then $E^+\cdot E^-=3$, and if the preimage of the node becomes intersection points of $E^+$ and $E^-$ then $E^+\cdot E^-=5$. Therefor we have $n=0$ (resp. 2) when $E^+\cdot E^-=3$ (resp. $E^+\cdot E^-=5$).
Finally, by Remark \ref{rem:restriction} (ii) and Lemma~\ref{lem:key-1} or Example 3.1 we have 
\[
\boldc(v_0, s(E^+))={}^t\left[\begin{array}{cccc} 0 & 0 & 0 & 0 \end{array}\right], \boldc(v_i, s(E^+))=0\,(i=1,2,3), 
\]
which implies that $P_\tau=0$ by the arguments in Section \ref{sec:torsion} . 

So far, we have two possibilities
\[
P_{E^+}=
\begin{cases}
O \\
2P_o
\end{cases}
\]

We  can determine which case occurs as follows:

Suppose $s(E^+)=O$. This implies that for every smooth singular fiber $F$, the  sum of the three intersection points of $E^+$ and $F$ under the group operation of $F$ must equal
 $O\cap F$. Let $L_{12}$ be the line through $q_1, q_2$ and choose $z_o$ to be the intersection of $L_{12}$ and $L_0$. Then the intersection points of $E^+$ and  the smooth fiber corresponding to $L_{12}$, lying above $q_1, q_2$ become 2-torsion points on $F$. This and  the fact that the sum of these two points with the third intersection point must equal $O\cap F$  implies that the third intersection point must also be a 2-torsion point on $F$, hence must lie over $q_3$. Hence $q_3\in  L_{12}$ and $q_1, q_2, q_3$ become collinear.

On the other hand, if $s(E ^+)=s(2P_o)$, we have $s(E^+)\cdot O=0$. This implies that for every smooth fiber $F$, the sum of the intersection points of $E^+$ and $F$ cannot equal $O$ under the group operation of $F$. Hence in this case, $q_1, q_2, q_3$ cannot be collinear.

Summing up, we have the following proposition.
\begin{prop}\label{prop:imageE} Let $E$ be a smooth or nodal cubic having the combinatorics given above. Then, the following statements hold:
\begin{enumerate}
\item If $E$ is a smooth cubic and is splitting then  $s(E^+)=O$ and $q_1, q_2, q_3$ are collinear.
\item If $E$ is a nodal cubic, then $E$ always splits and
\begin{enumerate}
\item $s(E^+)=O$ if and only if $q_1, q_2, q_3$ are collinear, and
\item  $s(E^+)=s_{(2P_o)}$ otherwise.
\end{enumerate}
\end{enumerate}
\end{prop}
\proof All the statements follow from the above discussions. We note that if $E$ is a nodal cubic, $E$ is a splitting curve as it is rational.

%
  \section{Dihedral covers and proof of Theorem~\ref{thm:refine}}

\subsection{Dihedral covers}

Let $D_{2n}$  be the dihedral group of order $2n$.
In order to prove Theorem~\ref{thm:refine}, we consider the existence/non-existence for Galois covers of $\PP^2$
whose Galois group is isomorphic to $D_{2n}$ ($D_{2n}$-covers, for short). In this
section, we summarize some facts on $D_{2n}$-covers which we need for our proof.
As for general terminologies on Galois covers, we refer to \cite[Section 3]{act}, or  \cite{tokunaga14}.
As for $D_{2n}$-covers, we refer \cite{artal-tokunaga, bannai-tokunaga15, bannai-tokunaga17, tokunaga94, tokunaga14}.
 In our previous papers such as   \cite{bannai-tokunaga15,  bannai-tokunaga17,tokunaga14},
we considered the cases when $n$ is odd. In this article, however, we also consider the case where $n$ is even. We here 
introduce some results on the even case based on  \cite[Proposition~0.6 and Remark~3.1]{tokunaga94} as follows:

\begin{prop}\label{prop:suff}{Let $n$ be an even integer $\ge 4$. Let $f :  S \to \Sigma$ be a smooth finite double cover of a simply connected smooth projective surface $\Sigma$. Let $\sigma_f$ be the 
involution on $S$ determined by the covering transformation for $f$. Let $C, D$ and $D_o$ be divisors on $S$ satisfying the following 
properties:

\begin{enumerate}

\item[(i)]  If  $C$ is irreducible such that $\sigma^*_f C \neq C$

\item[(ii)] $D$ is a reduced divisor or $D = \emptyset$. If $D \neq \emptyset$ or each irreducible component $D_j$, there exists an irreducible divisor $B_j$ on $\Sigma$ such that $D_j = f^*B_j$.

\item[(iii)]  $C + n/2D- \sigma_f^*C \sim nD_o$.

\end{enumerate}

Then there exists a $D_{2n}$-cover $\pi : X \to \Sigma$ such that $(a)$ $D(X/\Sigma) = S$ and 
$(b)$ $\Delta_{\pi} = \Delta_f \cup f(\Supp(C + D))$.
}
\end{prop}

\proof We rewrite $D_o$ as a difference of effective divisors: $D_o = D^+ - D^-$. Now put
$D_1 = C, D_2 = D, D_3 = D^-$ and $D_4 = D^+$. By \cite[Remark~3.1]{tokunaga94}, the condition (e) in \cite[Propositons~0.6]{tokunaga94} 
is satisfied. Hence, by Proposition~0.6, we have a $D_{2n}$-cover as desired. \qed

As for the ^^ ^^ converse" for Proposition~\ref{prop:suff}, it can be stated as follows:

\begin{prop}\label{prop:necc}{Let $f : S \to \Sigma$,  $\sigma_f$ and $D$ be as in Proposition~\ref{prop:suff}. If there exists a $D_{2n}$-cover
($n$: even $\ge 4$) $\pi : X \to \Sigma$ such that
\begin{enumerate}
\item[(i)] $D(X/\Sigma) = S$ and
\item[(ii)] $\pi$ is branched at $2(\Delta_f + f(D)) + n\overline{C}$.
\end{enumerate}
Then:
\begin{enumerate}
\item[(a)] $f^*\overline{C}$ is of the form $C + \sigma_f^*C$, $C \neq \sigma_f^*C$.
\item[(b)] There exists divisor $D_o$ on $S$ such that
$C + n/2 D - \sigma_f^*C \sim nD_o$.
\end{enumerate}
}
\end{prop}

\proof Let $D_1, D_2, D_3$ and $D_4$ be the divisors in \cite[Proposition~0.7]{tokunaga94}. Then the conditions (i) and (v) in
 \cite[Proposition~0.7]{tokunaga94} and (ii) as above, we may assume:
 \[
 D_1 = aC, \, 0 < a < (n-1)/2, \, D_2 = D, \, D'_o = D_4 - D_3,
 \]
 and we have $aC - n/2 D - a\sigma_f^*C \sim nD_o$. Take $u \in \CC(S)$ such that $(u) = aC - n/2 D - a\sigma_f^*C - nD_o$.
 By the proof of \cite[Proposition~0.7]{tokunaga94}, $X$ is the $\CC(S)(\sqrt[n]{u})$-normalization of $S$. This shows that
 the ramification index along $C$ is given by $n/\gcd(a, n)$, i.e., $\gcd(a, n) = 1$ by the condition (ii). Let $k$ be an integer such that
 $ak \equiv 1 \bmod n$. Then
 \begin{eqnarray*}
  & & akC + \frac{kn}2 D - ak\sigma_f^*C \\
  &= & C + \frac n2 D - \sigma_f^*C + n\left\{\frac {nk-1}n (C - \sigma_f^*C) + \frac {k-1}2 D\right \}
  \end{eqnarray*}
  Hence we have
  \[
  C + \frac n2 D - \sigma^*_fC \sim n\left\{D'_o - \frac {nk-1}n (C - \sigma_f^*C) - \frac {k-1}2 D\right\}.
  \]
 \qed

\subsection{Proof of Theorem~\ref{thm:refine}}

We first remark that, as  $S_{\mcQ, z_o}$ is
a rational elliptic surface, there is no difference between algebraic and linear equivalence.  For the covering 
transformation of $f_{\mcQ, z_o}$, we denote it by $\sigma_{\mcQ, z_o}$ for simplicity.

{\bf Proof for  (i).} In this case, $q_1, q_2$ and $q_3$ are collinear. By Proposition~\ref{prop:imageE}, we have
$s(E^+) = O$, i.e., 
\[
E^+\sim 3O+ 3F-2\Theta_{\infty,1}-2\Theta_{\infty,2}-2\Theta_{\infty,3}-3\Theta_{\infty,4}.
\]
$O, F$ and $\Theta_{\infty,i}$ are invariant under $\sigma_{\mcQ, z_o}$,  $E^+ \sim E^-$. Hence, by
\cite[Proposition~1.1, Corollary~1.2]{artal-tokunaga}, there exists an $n$-cyclic cover $g_n : \hat{X}_n \to S_{\mcQ, z_o}$ such that 
(i) $\Delta_{g_n} = E^+ \cup E^-$ and 
(ii) $f_{\mcQ, z_o}\circ g_n : \widehat{X}_n \to (\widehat{\PP^2})_{z_o}$ is a $D_{2n}$-cover. The stein  factorization  $X_n$ of $q\circ q_{z_o}\circ f_{\mcQ, z_o}\circ g_n :
\widehat{X}_n \to \PP^2$ is a $D_{2n}$-cover of $\PP^n$ as desired.

\medskip

{\bf Proof for (ii).}  We prove the existence of a $D_8$-cover.   Let us recall the notation: $L_o\cap L_i = \{p_i\}$ and $L_i\cap L_j = \{p_{ij}\}$.
Let $\overline{p_ip_{jk}}$ ($\{i, j, k\} = \{1, 2,3\}$) be the lines connecting $p_i$ and $p_{jk}$. Each $\overline{p_ip_{jk}}$ gives rise to sections $s_i^{\pm}$ such
that $\langle P_{s_i}^{\pm}, P_{s_i}^{\pm}\rangle = 1/2$. We choose one of them as a generator $s_o$ of the free part $A_1^*$ (note that the differences between
$s_i^{\pm}$ and $s_j^{\pm}$ are a translation-by-$2$-torsions). Choose $P_{s_1^+}$ as $P_{o}$. Note that we change $\pm$ so that $s(E^+) = s_{(2P_o)}$ if necessary. Then as
\[
2s_o \sim s_{(2P_o)} + O + F - 2\Theta_{\infty,1} - \Theta_{\infty,2} - \Theta_{\infty, 3} - 2\Theta_{\infty, 4} - \Theta_{1,1},
\]
we have 
\begin{eqnarray*}
E^+ & \sim & s_{(2P_o)} + 2O + 2F - 2\Theta_{\infty, 1} - 2\Theta_{\infty, 2} - 2\Theta_{\infty, 3} - 3\Theta_{\infty, 4} \\
       & \sim & 2s_o + O + F - \Theta_{\infty, 2} - \Theta_{\infty, 3} - \Theta_{\infty, 4} + \Theta_{1, 1}.
\end{eqnarray*}
Since $\sigma_{\mcQ, z_o} = [-1]_{\varphi_{\mcQ, z_o}}$,  we have  $s(E^-) = s_{(-2P_o)}$. Also as     
\[
-2s_o \sim s_{(-2P_o)} -3O -3F + 2\Theta_{\infty,1} + \Theta_{\infty,2} + \Theta_{\infty, 3} + 2\Theta_{\infty, 4} + \Theta_{1,1},
\]
we have
\[
E^-  \sim  -2s_o + 5(O + F) - 4\Theta_{\infty, 1} - 3\Theta_{\infty, 2} - 3\Theta_{\infty, 3} - 5\Theta_{\infty, 4} - \Theta_{1,1}.
\]
Thus we have
\[
E^+ + 2(\Theta_{\infty, 2} + \Theta_{\infty, 3} + \Theta_{1,1}) - E^- \sim 4(s_o - O - F + \Theta_{\infty, 1} + \Theta_{\infty, 2} + \Theta_{\infty, 3} + \Theta_{\infty,4} + \Theta_{1,1}).
\] 

Now by Proposition~\ref{prop:suff}, we have a $D_8$-cover $\widehat{\pi}_4: \widehat{X}_4 \to (\widehat{\PP^2})_{z_o}$ branched at $2\Delta_{f_{\mcQ, z_o}} + 4(f_{\mcQ, z_o}(E^+ \cup
\Theta_{\infty, 2} \cup\Theta_{\infty, 3} \cup \Theta_{1,1}))$ and the Stein factorization of $q\circ q_{z_o}\circ \hat{\pi}_4$ gives the desired $D_8$-cover
$\pi_4 : X_4 \to \PP^2$ branched at $2Q + 4E$.

\medskip
 
 We now go on to prove non-existence of $D_{2n}$-covers as in Theorem~\ref{thm:refine} for $n \ge 3$ except $n = 4$.
 Suppose that there exists a $D_{2n}$-cover $\pi_n : X_n \to \PP^2$ of $\PP^2$ as in Theorem~\ref{thm:refine}. Then
 we have $D(X_n/\PP^2) = S'_{\mcQ}$ and $f'_{\mcQ} = \beta_1(\pi_n)$ Let $\widehat{X}_n$ be the $\CC(X_n)$-normalization of 
 $(\widehat{\PP^2})_{z_o}$. Then we have the following:
 
 \begin{itemize}
 
 \item The induced morphism $\hat{\pi}_n : \widehat{X}_n \to (\widehat{\PP^2})_{z_o}$ is a $D_{2n}$-cover of 
 $(\widehat{\PP^2})_{z_o}$.
 
 \item  $D(\widehat{X}_n/(\widehat{\PP^2})_{z_o} = S_{\mcQ, z_o}$ as $S_{\mcQ, z_o}$ is the $\CC(S'_{\mcQ})$-normalization of
 $(\widehat{\PP^2})_{z_o}$.
 
 \item The image of the branch locus $\Delta_{\beta_2(\hat{\pi}_n)}$ of the form $f_{\mcQ, z_o}(E^+) + \overline {\Xi}$, where
 $\overline{\Xi}$ is contained in the exceptional set of $q\circ q_{z_o}$.
 
 \end{itemize}

\begin{lem}\label{lem:proof-of-refine}
\begin{enumerate}

 \item[(i)] If $n$ is odd, $\overline{\Xi} = \emptyset$.
 
 \item[(ii)] If $\overline{\Xi} \neq \emptyset$, then $n$ is even and the ramification index along $\Xi$ is $2$.
 
\end{enumerate}
\end{lem}

\proof   Let $E_o$ be an arbitrary irreducible component of $\overline{\Xi}$. In our case, $f^*_{\mcQ, z_o}(E_o)$ is
irreducible.
Now (i) follows from \cite[Corollary~2.4]{tokunaga14}, and (ii) from \cite[Proposition~0.7]{tokunaga94}. \proofend

\medskip

{\sl $n = $odd:} Choose any odd prime $p$ dividing $n$. Since we have a surjective morphism $D_{2n} \to D_{2p}$,
we have a $D_{2p}$-cover $\hat{\pi}_p: \widehat{X}_p \to (\widehat{\PP^2})_{z_o}$ such that 
$D(\widehat{X}_p/(\widehat{\PP^2})_{z_o}) = S_{\mcQ, z_o}$ and $\beta_2(\hat{\pi}_p)$ is branched 
at $p(E^+ + E^-)$ by Lemma~\ref{lem:proof-of-refine}. By
\cite[Theorem~3.2]{tokunaga14}, $s(E^+)$ is $p$-divisible in $\MW(S_{\mcQ, z_o})$,   but this contradicts to  Proposition~4.1.

{\sl $n =$ even:} By Propositon~5.2,  there exist divisors $D$ and $D_o$ on $S_{\mcQ, z_o}$ such that
\[
E^+ + \frac n2 D - E^- \sim nD_o,
\]
where $\Supp(D)$ is contained in the exceptional set of $q\circ q_o$ by Lemma~\ref{lem:proof-of-refine}.
By \cite[Theorem~1.3]{shioda90}, this implies $P_{E^+}-P_{E^-}$ is $n$-divisible in $\MW(S_{\mcQ, z_o})$.
On the other hand, by Proposition~4.1,  we have
\[
P_{E^+} - P_{E^-} = 2P_o - (-2P_o) = 4P_o.
\]
This proves the non-existence for $n \ge 6$.

\begin{rem}{\rm As we noticed before, there are $4$ different choices for $s_o$.  
If we choose  generator different from $s_o$,  the divisor $\Theta_{\infty, 2} + \Theta_{\infty, 3} + \Theta_{1,1}$ changes.
This means that we have $4$ different $D_8$-covers in the case of curves of Type II.

}
\end{rem}

%

\noindent Shinzo BANNAI\\
National Institute of Technology, Ibaraki College\\
866 Nakane, Hitachinaka-shi, Ibaraki-Ken 312-8508 JAPAN \\
{\tt sbannai@ge.ibaraki-ct.ac.jp}\\

\noindent Hiro-o TOKUNAGA\\
Department of Mathematics and Information Sciences\\
Tokyo Metropolitan University\\
1-1 Minami-Ohsawa, Hachiohji 192-0397 JAPAN \\
{\tt tokunaga@tmu.ac.jp}


\begin{thebibliography}{99}
   
   \bibitem{act} E.~Artal Bartolo, J.-I.~Codgolludo and H.~Tokunaga:
   \emph{A survey on Zariski pairs}, Adv.Stud.Pure Math., \textbf{50}(2008), 1-100.
   
   
  
  
  %
  %
%
   \bibitem{artal-tokunaga} E.~Artal Bartolo and H.~Tokunaga: \emph{Zariski $k$-plets of rational curve arrangements and dihedral covers},  Topology Appl. 
  \textbf{142} (2004),  227-233.
  
  %
  \bibitem{bannai-tokunaga15} S.~Bannai and H.~Tokunaga: \emph{Geometry of bisections of elliptic surfaces and Zariski N-plets for conic arrangements},  Geom. Dedicata {\bf 178}(2015), 219 - 237.
  
 \bibitem{bannai-tokunaga17} S.~Bannai and H.~Tokunaga: \emph{Geometry of bisections of elliptic surfaces and Zariski N-plets II},  Topology and its Applications, {\bf 231}(2017), 10 - 25
  
  %
%
%
  
  %
  
  \bibitem{bty17} S.~Bannai, H.~Tokunaga and M.~Yamamoto: Rational points of elliptic surfaces 
and 
Zariski $N$-ples for cubic-line, cubic-conic-line arrangements, arXiv:1710.02691
  
  
   
%
%


\bibitem{cogo-kloosterman} J.-I.~Cogolludo-Agusttin and R.~Kloosterman: \emph{Mordell-Weil groups and Zariski triples. Geometry and arithmetic}, 75-89, EMS Ser. Congr. Rep., Eur. Math. Soc., Z\"urich, 2012. 



\bibitem{cogo-lib}  J. -I.~Cogolludo-Agustin and A.~Libgober: \emph{Mordell-Weil groups of elliptic threefolds and the Alexander module of plane curves}, 
J. Reine Angew. Math. {\bf 697}(2014), 15-55.



 
 \bibitem{horikawa} E.~Horikawa: \emph{ On deformation of quintic surfaces},
\rm Invent. Math. {\bf 31} (1975), \rm $43 - 85$.


\bibitem{kloosterman2013} R.~Kloosterman: \emph{Cuspidal plane curves, syzygies and a bound on the MW-rank},  J. Algebra {\bf 375} (2013), 216-234

\bibitem{kloosterman2014} R.~Kloosterman: \emph{Mordell-Weil lattices and toric decompositions of plane curves},  Math. Ann. {\bf 367} (2017), 755-783.


\bibitem{kodaira} K.~Kodaira: \emph{On compact analytic surfaces II-III}, Ann. of Math. \textbf{77}
(1963), 563-626, \textbf{78}(1963), 1-40..

%

\bibitem{lib} A.~Libgober:  \emph{On Mordell-Weil groups of isotrivial abelian varieties over function fields},  Math. Ann. {\bf 357} (2013), 605-629.

 
 \bibitem{miranda-basic} R.~Miranda: \emph{Basic theory of elliptic surfaces}, Dottorato di Ricerca in Matematica, ETS Editrice, Pisa, 1989.

 
\bibitem{miranda-persson} R.~Miranda and U.~Persson: \emph{
On extremal rational elliptic surfaces}, Math. Z. \textbf{193}(1986), 537-558.

%
%
\bibitem{oguiso-shioda}  K.~Oguiso and T.~Shioda: \emph{The Mordell-Weil lattice of Rational
Elliptic surface}, Comment. Math. Univ. St. Pauli \textbf{40}(1991), 83-99.

%
%


 \bibitem{shioda90} T.~Shioda: \emph{On the Mordell-Weil lattices}, \rm Comment. Math. Univ. St. Pauli
\textbf{39} (1990), 211-240.  


%


\bibitem{tokunaga94} H.~Tokunaga:  \emph{On dihedral Galois coverings}, \rm Canadian J. of
Math. {\bf 46} \rm (1994),1299 - 1317.

\bibitem{tokunaga98}  H.~Tokunaga: \emph{Some examples of Zariski pairs arising from certain 
elliptic K3 surfaces}, \rm
Math. Z. {\bf 227} (1998), 465-477, \emph{II},  \rm Math.Z. {\bf 230} (1999), 389-400
%
%
%
\bibitem{tokunaga04} H.~Tokunaga: \emph{Dihedral covers and an elemetary arithmetic on elliptic surfaces}, 
J. Math. Kyoto Univ. \textbf{44}(2004), 55-270.


%
%
\bibitem{tokunaga12} H.~Tokunaga: \emph{Some sections on rational elliptic surfaces and certain special
conic-quartic configurations},  Kodai Math. J.  \textbf{35}(2012), 78-104.

\bibitem{tokunaga14} H.~Tokunaga: \emph{Sections of elliptic surfaces 
and Zariski pairs for conic-line arrangements via dihedral covers},  J. Math. Soc. Japan \textbf{66}(2014) 613-640.


%
%
%

%





%


%



%
%




%



%
%

%












%

%
%
%
%
%
%
%
%
%

\end{thebibliography}
  \end{document}